\documentclass[11pt]{article}
\begin{document}
\title{On $p$-adic intermediate Jacobians}
\author{Wayne Raskind\thanks{Partially supported by NSF grant 0070850,
SFB 478 (M{\"u}nster), CNRS France, and sabbatical leave from the
University of Southern California}\hspace{.1in}
\setcounter{footnote}{6} Xavier Xarles\thanks{Partially supported
by grant BFM2003-06092 from DGI}}

\date{}
\maketitle

\newcommand{\free}[1]{{#1}/tors}

\newcommand{\bC}{{\bf C}}
\newcommand{\bZ}{{\bf Z}}
\newcommand{\bQl}{{\bf Q}_{\ell}}
\newcommand{\bZl}{{\bf Z}_{\ell}}
\newcommand{\Xb}{{\overline X}}
\newcommand{\Zb}{{\overline Z}}
\newcommand{\Eb}{{\overline E}}
\newcommand{\Kb}{{\overline K}}
\newcommand{\Cb}{{\overline C}}

\newcommand{\cX}{{\cal X}}
\newcommand{\Yb}{{\overline Y}}
\newcommand{\bQ}{{\bf Q}}
\newcommand{\bQp}{{\bf Q}_p}
\newcommand{\bZp}{{\bf Z}_p}
\newcommand{\bCp}{{\bf C}_p}
\newcommand{\bR}{{\bf R}}
\newcommand{\cE}{{\cal E}}
\newtheorem{proposition}{Proposition}
\newtheorem{theorem}[proposition]{Theorem}
\newtheorem{lemma}[proposition]{Lemma}
\newtheorem{definition}[proposition]{Definition}
\newtheorem{remark}[proposition]{Remark}
\newtheorem{example}[proposition]{Example}
\newtheorem{corollary}[proposition]{Corollary}
\newtheorem{conjecture}[proposition]{Conjecture}
\newcommand{\plim}{\displaystyle{\lim_{\stackrel{\longleftarrow}{n}}}\,}
\newcommand{\bF}{{\overline F}}
\newcommand{\kb}{{\overline k}}
\newcommand{\Fb}{{\overline F}}
\newcommand{\cO}{{\cal O}}
\newcommand{\im}{\mbox{Im}}
\newcommand{\bG}{{\bf G}}

\abstract{For an algebraic variety $X$ of dimension $d$ with
totally degenerate reduction over a $p$-adic field (definition
recalled below) and an integer $i$ with $1\leq i\leq d$, we define
a rigid analytic torus $J^i(X)$ together with an Abel-Jacobi
mapping to it from the Chow group $CH^i(X)_{hom}$ of codimension
$i$ algebraic cycles that are homologically equivalent to zero
modulo rational equivalence. These tori are analogous to those
defined by Griffiths using Hodge theory over $\bC$.  We compare
and contrast the complex and $p$-adic theories. Finally, we
examine a special case of a $p$-adic analogue of the Generalized
Hodge
Conjecture.}\\

AMS classification:  14K30 (primary), 14K99 and 14F20 (secondary)
\parindent=0cm
\section*{Introduction}

Let $E$ be an elliptic curve over the complex number field $\bC$.
The group $E(\bC)$ of $\bC$-points of $E$ may be written as

$$E(\bC)=\bC/\Lambda,$$

where $\Lambda$ is a lattice of periods of the form $\bZ\oplus
\bZ\tau$, with $\tau$ in the upper half plane.  Setting
$q=e^{2\pi\,i\tau}$, we can also write:

$$E(\bC)=\bC^*/q^{\bZ}.$$
In 1959, Tate showed that one can make sense of this latter
description in a general complete valued field $K$ if $q$ is of
absolute value less than 1.  Such elliptic curves are now called
{\it Tate elliptic curves}, and if $K$ is nonarchimedian then they
have split multiplicative reduction modulo the maximal ideal of
the valuation ring of $K$. Let $\ell$ be a prime number,
$E[\ell^n]$ be the group of points of order $\ell^n$ on $E$ over
an algebraic closure of $K$ and $T_{\ell}$ denote the Tate module

$$\plim E[\ell^n].$$

Then we have an exact sequence of $G$-modules, where $G$ is the
absolute Galois group of $K$:

$$0\to \bZl(1)\to T_{\ell}\to \bZl\to 0.$$

The class of this extension in $Ext^1_G(\bZl,\bZl(1))$ is given by
the class of $q$ in the completion of $K^*$ with respect to its
subgroups $K^{*\,\ell^n}$. The Galois module structure of such
extensions is reasonably simple but still very interesting.
Mumford, Raynaud and many others have generalized this theory to
other types of varieties such as abelian varieties with totally
multiplicative reduction and certain (Mumford) curves of higher
genus (see [Mu1,2] and [Ray1]). Manin and Drinfeld have given an
analytic definition of the Abel-Jacobi mapping in the case of Mumford curves [MD].\\

Recall that the important ingredients in the theory of $p$-adic
uniformization of an abelian variety $A$ (see for example [Ray1])
are two isogenous lattices of periods $\Lambda$ and $\Lambda'$ and
a nondegenerate pairing (where
$\Lambda'^{\vee}=\mbox{Hom}(\Lambda',\bZ))$:

$$\Lambda\times \Lambda'^{\vee}\to K^*.$$

Then the pairing gives a uniformization of $A$:

$$A\cong\mbox{Hom}(\Lambda'^{\vee},\bG_m)/\Lambda$$

as rigid analytic varieties.\\

Let now $K$ be a finite extension of $\bQp$ and $X$ be a smooth,
projective, geometrically connected variety over $K$.  Let $\Kb$
be an algebraic closure of $K$ and $\Xb=X\times_K\Kb$. In this
paper, we use the methods of our earlier paper [RX] to reinterpret
these lattices and the pairing in a purely algebro-geometric way
using comparison theorems between $p$-adic \'etale cohomology and
log-crystalline cohomology (that is, Tsuji's proof of the
semi-stable conjecture $C_{st}$ [Tsu], see also Faltings approach
in [Fa2]). The utility of these methods is that they generalize
quite easily to higher cohomology groups (with the lattices for
abelian varieties being the case of $H^1$), whereas the
uniformization methods have met
with somewhat less success here (but see [dS]).\\

In [RX], we formalized the notion of a variety $X$ over $K$ with
{\it totally degenerate reduction} (we will recall the definition
below) and we showed that for {\it all} prime numbers $\ell$, the
 \'etale cohomology groups $H^*(\Xb,\bQl)$  are, after a finite
 unramified extension,
successive extensions of direct sums of Galois modules of the form
$\bQl(r)$ for various $r$.  The main result of the present paper
associates to each odd \'etale cohomology group
$\prod_{\ell}H^{2i-1}(\Xb,\bQl)$ of such a variety a $p$-adic
analytic torus $J^i(X)$ of dimension equal to the Hodge number
$h^{i,i-1}$ and an Abel-Jacobi mapping:

$$CH^i(X)_{hom}\to J^i(X)(K)\otimes_{\bZ}{\bQ},$$
where the group on the left is that of cycles of codimension $i$
that are homologically equivalent to zero, modulo rational
equivalence (see Remark 17(i) in \S 3 below for why we have to
tensor with $\bQ$ on the right). We have that $J^1(X)$ is
 the Picard variety, $Pic^0(X)$, and $J^d(X)$  is the Albanese variety,
$Alb(X)$, where $d$ is the dimension of $X$.\\

We call the $J^i(X)$ ``$p$-adic intermediate Jacobians'', and we
expect them to provide a useful first step in the study of
algebraic cycles on such varieties, in a similar way as the
intermediate Jacobians of Griffiths [Gr] do for varieties $X/\bC$.
Our original motivation for studying these $p$-adic intermediate
Jacobians was to provide, at least in some cases, an
algebro-geometric interpretation for the intermediate Jacobians
$\tilde{J}^i(X)$ of Griffiths and the Abel-Jacobi map

$$CH^i(X)_{hom}\to \tilde{J}^i(X).$$

 These are defined by Griffiths in a complex analytic way using
the Hodge decomposition, and they have no known purely
algebro-geometric interpretation, in general.  Weil had earlier
introduced similar objects, but with a different complex structure
[W1,2]
\\

Recall that $\tilde{J}^i(X)$ is a (compact) complex torus, not
necessarily an abelian variety, whose  dimension is one half of
$B_{2i-1}=\mbox{dim}_{\bC}H^{2i-1}(X,\bC)$.  For $X$  a generic
quintic 3-fold, Griffiths (loc. cit.) used the Abel-Jacobi map to
$\tilde{J}^2(X)$ to detect codimension two cycles that are
homologically equivalent to zero, but not algebraically equivalent
to zero.\\

The big advantage of our $J^i$ is that it is defined over $K$ and
therefore has arithmetic meaning and applications.  A disadvantage
is that the dimension of $J^i$ is in general strictly smaller than
that of the analogous complex intermediate Jacobian $\tilde{J^i}$.
To help compensate for this deficiency, we relate the ``missing
pieces'' to higher odd algebraic K-theory, which we are not able
to do at present in the complex case (see Remark 17 (ii) for more
on this). In fact, the methods of this paper suggest going back to
the complex case and trying to reinterpret $\tilde{J^i}$ in a
similar way, using the idea that complex varieties should have
(totally degenerate) ``semi-stable'' reduction in the same way as
the varieties we study here do over complete discretely valued
fields (see [Ma] for what this means in the case of curves). In
work in preparation, the second author and Infante compare and
contrast the $p$-adic and complex cases for abelian varieties. In
[R2], the first author formulates a ``generalized Hodge-Tate
conjecture'' for varieties of the type considered in this paper
and in [R3], he proves a slightly weaker form of this conjecture
for divisors. We consider a special case of the generalized
Hodge-Tate conjecture in \S 4 below.
\\

This paper is organized as follows: in \S 1, we recall the
necessary preliminaries from [RX].  We then define the
intermediate Jacobians in \S 2 and the Abel-Jacobi mapping in \S
3.  Finally, in \S 4 we consider an example of the product of two
and then three Tate elliptic curves, and we determine in some
cases the dimension of the image in $J^2(X)$ of the Abel-Jacobi
mapping restricted to cycles algebraically equivalent to zero.
This dimension conjecturally depends on the rank of the space of
multiplicative relations between the Tate parameters of the
curves.   This is a special case of the generalized Hodge-Tate
conjecture mentioned above.  Our analysis shows that $J^2(X)$ is
never equal to this image. This agrees with the complex case,
since there we have that this image is contained in
$H^{1,2}/[H^{1,2}\cap H^3(X,{\bZ})]$, and is all of this last
group iff the three elliptic curves are isogenous and have complex
multiplication (see [Li], p. 1197). This cannot happen for Tate
elliptic curves, since CM curves have everywhere potentially good
reduction. \\

        The authors would like to thank, respectively, the Universitat
Aut{\`o}noma de Barcelona and the University of Southern
California for their hospitality.  This paper was completed while
the first author enjoyed the hospitality of Universit{\'e} de
Paris-Sud, SFB Heidelberg and Universit{\'e} Louis Pasteur
(Strasbourg). We also thank J.-L. Colliot-Th{\'e}l{\`e}ne and A.
Werner for helpful comments and information.

 \section{Notation and Preliminaries}
 \subsection{Notation}
         Let $K$ be a finite extension of $\bQp$, with valuation ring $R$
and residue field $F$, and let $X$ be a smooth, projective,
geometrically connected variety over $K$.  We denote by $\Kb$ an
algebraic closure of $K$, $\Xb=X\times_K\Kb$ and $G=Gal(\Kb/K)$.
If $r$ is a nonnegative integer, then $\bZl(r)$ denotes the Galois
module $\bZl$, Tate-twisted $r$ times.  If $r<0$, then
$\bZl(r)=\mbox{Hom}(\bZl(-r),\bZl)$.  $\bQl(r)$ is defined by tensoring these by $\bQl$.\\

By an {\it analytic torus} over $K$, we mean the quotient of an
algebraic torus over $K$ by a lattice as rigid analytic varieties (see for example [Ray1]
for the exact meaning of this).\\

Let $S$ be any domain with field of fractions $\mbox{frac}(S)$ of
characteristic zero. If $M$ is an $S$-module, we denote by
$\free{M}$ the torsion free quotient of $M$.  That is,
$M/M_{tors}$, where $M_{tors}$ denotes the submodule of torsion
elements in $M$. We say that a map between $S$-modules is an
isomorphism modulo torsion if it induces an isomorphism between
the torsion free quotients. If $M$ and $M'$ are torsion free
$S$-modules of finite rank, we say that a map $\phi\colon M\to M'$
is an isogeny if it is injective with cokernel of finite exponent
(as abelian group). In this case, there exists a unique map
$\psi\colon M'\to M$, the dual isogeny, such that $\phi\circ
\psi=[e]$ and $\psi\circ \phi=[e]$, where $e$ is the exponent of
the cokernel of $\phi$ and $[e]$ denotes the map multiplication by
$e$. In general, if $M$ and $M'$ are $S$-modules with torsion free
quotients of finite rank, we say that a morphism $\psi\colon M'\to
M$ is an isogeny if the induced map on the
torsion free quotients is an isogeny.\\

\subsection{Totally Degenerate Reduction}
Let $X$ be a smooth projective geometrically connected variety
over $K$. We assume that $X$ has a regular proper model $\cX$ over
$R$ which is strictly semi-stable, which means that the following conditions hold:\\

(*)   Let $Y$ be the special fibre of $\cX$.  Then $Y$ is reduced;
write
 $$Y=\bigcup_{i=1}^{n}Y_i,$$
 with each $Y_i$ irreducible. For each nonempty subset
$I=\{i_1,\dots,i_k\}$ of $\{1,\dots,n\}$, we
 set
 $$Y_I=Y_{i_1}\cap...\cap Y_{i_k},$$
scheme theoretically.  Then $Y_I$ is smooth and reduced over $F$
of pure codimension $|I|$ in $\cX$ if it is nonempty. See [dJ]
2.16 and [Ku], \S 1.9, 1.10, for a clear summary of these
conditions, as
well as comparison with other notions of semi-stability.\\

Let $\bF$ be an algebraic closure of $F$. We set
$\Yb_I=Y_I\times_F\bF$.  Note that these need not be connected.
\begin{definition}
We say that $Y$ is {\it totally degenerate over $F$} if it is
projective and the following conditions are satisfied for each
$Y_I$:

 \begin{itemize}
 \item[a)] For every $i=0,...,d,$ the Chow groups $CH^i(\Yb_I)$ are
 finitely generated abelian groups.

The groups $CH^i_{\bQ}(\Yb_I)$ satisfy the Hodge index
 theorem:  let
$$\xi:  CH^i(\Yb_I)\to CH^{i+1}(\Yb_I)$$
be the map given by intersecting with the class of a hyperplane
section for a fixed embedding of $Y$ in a projective space.
Suppose  we have $x\in CH^i_{\bQ}(\Yb_I)$ such that $\xi
^{d-i}(x)=0$. Then we have $(-1)^i\mbox{tr}(x\xi^{d-2i}(x))\geq
0$, with equality iff $x=0$.

 \item[b)] For every prime number $\ell$ different from $p$, the
{\'e}tale cohomology groups $H^{2i+1}(\Yb_I, \bZl)$ are torsion,
and the cycle map induces an isomorphism
$$\free{CH^i(\Yb_I)}\otimes \bZl \cong \free{H^{2i}(\Yb_I, \bZl (i))}.$$

Note that this is compatible with the action of the Galois group.

\item[c)] If $p>0$, let $W$ be the ring of Witt vectors of $\bF$.
Denote by $H^*(\Yb_I/W)$ the crystalline cohomology groups of
$\Yb_I$. Then the groups $H^{2i+1}_{\rm crys}(\Yb_I/W)$ are
torsion, and
 $\free{CH^i(\Yb_I)\otimes W(-i) }\cong \free{H^{2i}_{\rm crys}(\Yb_I/W)}$
via the cycle map.  Here $W(-i)$ is $W$ with the action of
Frobenius multiplied by $p^{i}$.

\item[d)]  $Y$ is ordinary, in that $H^r(\Yb,B\omega^s)=0$ for all
$i, r$ and $s$. Here $B\omega$ is the subcomplex of exact forms in
the logarithmic de Rham complex on $\Yb$ (see e.g. [I3],
D{\'e}finition 1.4).  By [H] and ([I2], Proposition 1.10, this is
implied by the $\Yb_I$ being ordinary in the usual sense, in that
$H^r(\Yb_I,d\Omega^s)=0$ for all $I, r$ and $s$.  For more on the
condition of ordinary, see [I2],
 Appendice and [BK1], Proposition 7.3)
\end{itemize}
\end{definition}

We will say that $X$ has \textit{totally degenerate reduction} if
it has a regular proper model $\cX$ over $R$ which is strictly
semi-stable and whose special fibre $Y$ is totally degenerate over
$F$. If $Y$ is totally degenerate and the natural maps

$$CH^i(Y_I)\to CH^i(\Yb_I)$$

are isomorphisms modulo torsion for all $I$, we shall say that $Y$
is \textit{ split totally degenerate}.  Since the $CH^i(\Yb_I)$
are all finitely generated abelian groups and there is a finite
number of them, there is a finite extension of the field of
definition where all the cycles given basis for any of the
$CH^i(\Yb_I)$ modulo torsion are defined. So, after a finite
extension, any totally degenerate variety becomes split totally
degenerate.
\\\\

  Examples of varieties with totally degenerate reduction
include abelian varieties $A$ such that the special fibre of the
N{\'e}ron model of $A$ over the ring of integers $R$ of $K$ is a
torus, and products of Mumford curves or other $p$-adically
uniformizable varieties, such as Drinfeld modular varieties [Mus]
and some unitary Shimura varieties (see Example 1  in [RX] for
further details and references).  We
expect that condition c) implies b).  \\

\subsection{Chow complexes}
For more details on this section, please see ([RX], \S 3).  We
write $Y^{(m)}:=\bigsqcup Y_I$, where the disjoint union is taken
over all subsets $I$ of $\{1,...,n\}$ with $\# I =m$ and $Y_I\ne
\emptyset$.

 For each pair of integers $(i,j)$ we define
 $$C^i_j(Y)=\bigoplus_{k\geq\max\{0,i\}} CH^{i+j-k}(Y^{(2k-i+1)}).$$
Note that there are only a finite number of summands here, because
$k$ runs from $\max\{i,0\}$ to $i+j$. Note also that $C^i_j(Y)$
can be non-zero only if $i=-d,...,d$ and $j=-i,...,d-i$. \\

Then define differentials $d^i_j:\, C^i_j(Y)\to C^{i+1}_j(Y)$ by
using the natural restriction and Gysin maps, and define
$T^i_j(Y)=\mbox{Ker}d^i_j/\mbox{Im}d^{i-1}_j$, the homology in
degree $i$ of the complex $C^i_j(Y)$. \\

The monodromy operator $N:\,C^i_j(Y)\to C^{i+2}_{j-1}(Y)$ is
defined as the identity map on the summands in common, and the
zero map on different summands. $N$ commutes with the
differentials, and so induces an operator on the $T^i_j(Y)$, which
we also denote by $N$. We have $N^i$ is the identity on $C^i_j(Y)$
for $i\geq 0$. The following result is a direct consequence of a
result of Guillen and Navarro ([GN] Prop. 2.9 and Th{\'e}or{\`e}me
5.2 or [BGS], Lemma 1.5 and Theorem 2), using the fact that the
Chow groups of the components of our $Y$ satisfy the hard
Lefschetz theorem and the Hodge index theorem.  It is a crucial
result for this paper.

\begin{proposition} Suppose that $Y$ satisfies the assumptions of \S 1.
Then the monodromy operator $N^i$ induces an isogeny:
$$N^i:\,T^{-i}_{j+i}(\Yb)\to T^{i}_{j}(\Yb)$$
for all $i\geq 0$ and $j$.
\end{proposition}

In the next theorem we summarize the main results in our first
paper on the \'etale cohomology of $X$ (see [RX], \S 4 Corollary 1
and \S 6, Theorem 3).

\begin{theorem} Let $X$ be a variety over $K$ of dimension $d$ with
totally degenerate reduction and with special fiber $Y$. Then
there is a monodromy filtration $M_{\bullet}$ on $H^*(\Xb,\bZl)$
for {\it all} $\ell$ such that we have a canonical isogeny
$$T^i_j(\Yb)\otimes\bZl(-j)\to
\mbox{Gr}^M_{-i}H^{i+2j}(\Xb,\bZl)$$ as $G$-modules, where $G$ is
the absolute Galois group of $K$. These isogenies are isomorphisms
for almost all $\ell$. We have also
$\mbox{Gr}^M_{-i+1}H^{i+2j}(\Xb,\bZl)$ is torsion. Moreover, for
$\ell \ne p$, this filtration coincides with the usual monodromy
filtration on $H^i(\Xb,\bQl)$,  and the map $N$ tensored with
$\bQl$ is the $\ell$-adic  monodromy map on the graded quotients.
For $\ell=p$, this filtration is obtained from the monodromy
filtration on the log-crystalline cohomology of $\Yb$ by applying
the functor $D_{st}$ (see ibid for more details).
\end{theorem}

\begin{remark} The idea of looking at the Chow
complexes was suggested by analogies with mixed Hodge theory, as
described in the work of Deligne [De1], Rapoport-Zink [RZ],
Steenbrink [St], and in work of Bloch-Gillet-Soul{\'e} [BGS] and
Consani [Con] on the monodromy filtration and Euler factors of
L-functions.
\end{remark}

\section{$p$-adic Intermediate Jacobians}

The goal of this section is to define, for every $j=1, \dots
d:=\dim(X)$, rigid analytic tori $J^j(X)$ associated to $\cX$ with
totally degenerate special fibre in the sense of \S 1.2. To do
this, we will look for a subquotient of the {\'e}tale cohomology
group $H^{2j-1}(\Xb,\bZl(j))$ that is an extension of $\bZl$'s by
$\bZl(1)$'s, as the $\ell$-Tate module of an analytic torus must
be. The monodromy filtration gives a natural way of finding such
subquotients.     Using Theorem 3, we see that the graded
quotients for this filtration have natural $\bZ$-structures.  The
torsion free quotients of these $\bZ$-structures will give us
lattices $\Lambda$ and $\Lambda'$ (to be defined below), and we
will define a nondegenerate period pairing
$$\Lambda\times \Lambda^{'\vee} \to K^*,$$
where $\Lambda^{'\vee}=\mbox{Hom}(\Lambda',\bZ)$.\\

Then $J^j(X)$ is defined to be
$\mbox{Hom}(\Lambda^{\vee},K^*)/\Lambda'$. To define the
Abel-Jacobi mapping, we use the groups $H_g^1(K,-)$ of Bloch-Kato
([BK2], \S 3), and the fact that
$$J^j(X)\cong \prod_{\ell}H^1_g(K,M_{1,\ell}/M_{-3,\ell}),$$
where $M_{\bullet,\ell}$ is the monodromy filtration on
$H^{2j-1}(\Xb,\bZl(i))$ modulo torsion (see Theorem 3 for this
filtration), and the product is taken over all prime numbers
$\ell$ . Then the Abel-Jacobi mapping is obtained from the
$\ell$-adic Abel-Jacobi mapping ([Bl], [R1]). It would be helpful
to define the Abel-Jacobi mapping by analytic means, so that it
could be defined directly over a field such as $\bCp$ (the
completion of an algebraic closure of $K$). For example, Besser
[Be] has interpreted the $p$-adic Abel-Jacobi mapping for
varieties with good reduction over $p$-adic fields in terms of
$p$-adic integration.\\

 Note also that our $J^j(X)$ depend {\it a
priori} on the choice of regular proper model $\cX$ of $X$ over
$R$. We expect that they can be defined directly from the
monodromy filtration on $X$, but we have not managed to do this.
We can show that modulo isogeny, they only depend on $X$
(see Remark 10(ii) in \S 2).\\

To shorten notation, we will denote $T^i_j(\Yb)$ by $T^i_j$. Since
all the $\free{T^{i}_{j}}$ are abelian groups of finite rank, the
action of the Galois group, which by definition of the $T^i_j$ is
unramified, will factor through a finite group. Suppose first of
all that the absolute Galois group of $F$ acts trivially on
$\free{T^{-1}_{j}}$ and on $\free{T^{1}_{j-1}}$ for any $j$; this
will happen after a finite unra\-mi\-fied extension of $K$. We
define a pairing for every $j$:
$$\{-,-\}_l:\free{T^{-1}_{j}}\times \free{T^{1\vee}_{j-1}}\to K^{*(\ell)}
:=\plim K^*/K^{*\ell^n}.$$  This is done by using properties of
the monodromy filtration $M_i(H)$ on $H:=H^{2j-1}(\Xb,\bZl(j))$.
The idea is to look at the quotient
$$\widetilde{M}:=\free{(M_1(H)/M_{-3}(H))},$$ which is a subquotient of
$H^{2j-1}(\Xb,\bZl(j))$ which is, modulo torsion, an extension of
$\bZl$'s by $\bZl(1)$'s. But due to the possible torsion we have
to do a slight modification of $\widetilde{M}$.\\

Consider the natural epimorphism from $\widetilde{M}$ to
$W_{1}:=\free{\mbox{Gr}_{1}^M H}$. Define then $W_{-1}$ as the
kernel of this map. Observe that $W_{-1}$ is isogenous to
$\free{\mbox{Gr}_{-1}^M H}$, and the na\-tu\-ral isogeny
$\psi\colon\free{\mbox{Gr}_{-1}^M H}\to W_{-1}$ has cokernel
isomorphic to the cokernel of $M_1(H)/M_{-3}(H)\!\to\!
\mbox{Gr}_{1}^M H=M_1(H)/M_{-1}(H)$ and with exponent $e_{\ell}$.
Considering the dual isogeny $\varphi$ from $W_{-1}$ to
$\free{\mbox{Gr}_{-1}^M H}$ such that $\psi \circ
\varphi=[e_{\ell}]$ (see \S 1.1), we can push out the extension
$$\: 0\to W_{-1}\to \widetilde{M} \to W_{1}\to 0,$$
with respect to this isogeny. Recalling that we have canonical
isogenies
$$\psi_1\colon \free{\mbox{Gr}_{-1}^M H}\to
\free{T^{1}_{j-1}}\otimes\bZl(1) \mbox{ and }\psi_2\colon
\free{\mbox{Gr}_{1}^M H}\to \free{T^{-1}_{j}}\otimes\bZl,$$ we
push out with respect to $\psi_1$ and pull back with respect to
the dual isogeny associated to $\psi_2$ to get an extension:
$$0\to \free{T^{1}_{j-1}}\otimes\bZl(1)\to E'_{\ell} \to
\free{T^{-1}_{j}}\otimes\bZl\to 0.$$ Observe that for almost all
$\ell$ the isogenies $\psi_1$ and $\psi_2$ are all isomorphism.

Denote by $e:=\prod_{\ell} e_{\ell}$; it is well defined because
$e_{\ell}=1$ for almost all $\ell$ since the graded quotients for
the monodromy filtration are torsion-free for almost all $\ell$
(this follows from Lemma 1 of [RX]). We use these numbers to
normalize these extensions: for any $\ell$ consider the extension
$E_{\ell}$ obtained from $E'_{\ell}$ by pushing out with respect
to the map
$$[e/e_{\ell}]:
\free{T^{1}_{j-1}}\otimes\bZl(1) \to
\free{T^{1}_{j-1}}\otimes\bZl(1).$$

\begin{lemma} Consider the $\bQl$-vector space $E_{\ell}\otimes_{\bZl}\bQl$
constructed above. Then, for $\ell\ne p$, we have that the
monodromy map is given by the composition
$$E_{\ell}\otimes_{\bZl}\bQl \to
T^{-1}_{j}\otimes\bQl\stackrel{\widetilde{N}}{\to}
T^{1}_{j-1}\otimes\bQl\to E_{\ell}\otimes_{\bZl}\bQl$$ where the
left and right hand side maps are the natural ones and the map
$\widetilde{N}$ is the composition of the map $N$ defined in \S
1.3 and the multiplication-by-$e$ map. For $\ell=p$, applying the
functor $D_{st}$ (see the introduction to [Tsu]) to this
composition we get the monodromy map on the corresponding filtered
$(\Phi,N)$-module.
\end{lemma}

{\bf Proof}: Consider first the case $\ell\ne p$. Observe that for
a $\bQl$-module $V$ with action of the absolute Galois group such
that monodromy filtration satisfies that $M_{i-1}=0 \subseteq
M_i\subseteq M_{i+1}\subseteq M_{i+2}=V$, then the monodromy map
$N$ is the composition of
$$V \to V/M_i \stackrel{N'}{\to} M_i(1) \to V,$$
where $N'\colon V/M_i \to M_i$ is the induced map by $N$.

Now by ([RX],\S 4, Corollary 2) if $X$ has totally degenerate
reduction, and $M_{\bullet}$ is the monodromy filtration on the
$\ell$-adic cohomology $H^{i+2j}(\Xb,\bQl)$, the induced map
$$\mbox{Gr}^M_{-i}H^{i+2j}(\Xb,\bQl) \stackrel{N'}\rightarrow
\mbox{Gr}^M_{-i-2}H^{i+2j}(\Xb,\bQl)(1)$$ coincides with the map
$N\colon T^{i}_{j}\otimes\bQl(-j)\to T^{i+2}_{j-1}\otimes\bQl(-j)$
that we defined in section \S2, by composing with the isomorphism
$$T^i_j(\Yb)\otimes\bQl(-j)\cong
\mbox{Gr}^M_{-i}H^{i+2j}(\Xb,\bQl)$$ deduced from the weight
spectral sequence and the cycle maps on the components of the
special fiber $\Yb^{(k)}$. In particular we have that, using the
notation in the construction above, the composition
$$\widetilde{M}\otimes_{\bZl}\bQl \to W_1\otimes_{\bZl}\bQl\cong T^{-1}_j\otimes\bQl
\stackrel{N}{\to} T^1_{j-1}\otimes\bQl\cong
W_{-1}\otimes_{\bZl}\bQl\to \widetilde{M}\otimes_{\bZl}\bQl$$ is
the monodromy map on $\widetilde{M}$.

To get the result observe that the $\bQl$-vector space
$E\otimes_{\bZl}\bQl$ was obtained from
$\widetilde{M}\otimes_{\bZl}\bQl$ by pushing out with respect to
the map multiplication by $e$ on the
$T^1_{j-1}\otimes_{\bZl}\bQl$, so we have a commutative diagram
$$\matrix{0\to& T^{1}_{j-1}\otimes\bQl(1)&\to& E_{\ell}\otimes_{\bZl}\bQl
&\to &T^{-1}_{j}\otimes\bQl&\to 0\cr & \uparrow [e]&&\uparrow g
&&\|&\cr 0\to& T^{1}_{j-1}\otimes\bQl(1)&\to&
\widetilde{M}\otimes_{\bZl}\bQl &\to &T^{-1}_{j}\otimes\bQl&\to
0,}$$ and the monodromy operator $N$ on $E_{\ell}$ is obtained
from the monodromy operator on $\widetilde{M}$ by composing with
the isomorphism $g$ and its inverse:
$$E_{\ell}\stackrel{g^{-1}}{\to}\widetilde{M}\stackrel{N}{\to}\widetilde{M}(-1)\stackrel{g}{\to}E_{\ell}(-1).$$

Finally, to obtain the result in the case $\ell=p$, one uses a
similar argument for the log-crystalline cohomology by applying
([RX], Corollary 3 and Theorem 3). \\

 Now let $\Lambda=\free{T^{-1}_j}$ and $\Lambda'=\free{T_{j-1}^1}$.  We are going to define a pairing
$$\Lambda\times \Lambda^{'\vee}\to K^*$$
  as follows:  let $\alpha\in\Lambda, \beta\in \Lambda^{'\vee}$, and  tensor the subgroups these
elements generate with $\bZl$ and $\bZl(1)$, respectively. Pull
back the extension $E_{\ell}$ (which was introduced just before
Lemma 5) by the former and push it out by the latter to get an
extension of $\bZl$ by $\bZl(1)$, hence an element of
$$\mbox{Ext}^1_G(\bZl,\bZl(1))=K^{*({\ell})},$$
the $\ell$-completion of the multiplicative group of $K$.  We
denote this element by $\{\alpha,\beta\}_{\ell}$ and the product
of the pairings over all $\ell$ by
 $\{\alpha,\beta\}$.  This pairing \textit{a priori} takes values
 in $\widehat{K^*}=\plim K^*/K^{*\,n}$.

\begin{proposition} The image of the pairing $\{-,-\}$ in $\widehat{K^{*}}$
is contained in the image of the natural map $K^* \to
\widehat{K^*}$.
\end{proposition}

In the proof, we shall use the following well-known lemma:
\begin{lemma}  Let $K$ be a finite extension of $\bQp$ with normalized
 valuation $v$.
  Denote by $\widehat{v}$ the ``completed valuation''
$$\widehat{K^*}\to \widehat{\bZ}.$$

Then the completed valuation of an element $\alpha\in
\widehat{K^*}$ lies in
 the image of $\bZ$ in $\widehat{\bZ}$ iff $\alpha$ lies in the image of
$K^*$ in $\widehat{K^*}$.
\end{lemma}

{\bf Proof of lemma}:  Let $R$ be the ring of integers in $K$.
Then since $K$ is finite over $\bQp$, the groups $R^*/R^{*n}$
 are finite for every integer $n$, and $R^*$ is complete with respect
 to the topology induced by its subgroups of finite index.  Thus the
natural map $R^*\to \widehat{R^*}$ is an isomorphism and the map
$K^*\to \widehat{K^*}$ is injective.  Hence the induced map

$$\widehat{K^*}/K^*\to \widehat{\bZ}/\bZ$$

is an isomorphism, which proves the lemma.

\begin{remark}
    The use of this lemma is the one place where our construction will
 not work over a larger field such as $\bCp$.  What seems to be needed is to
 define our extensions in some sort of category ${\cal C}$ of ``analytic
mixed motives''
 over $\bCp$, in which we would have $Ext^1_{\cal C}(\bZ,\bZ(1))=\bCp^*$.
As is well-known,
 extensions of $\bCp$ by $\bCp(1)$ in the category of continuous
$Gal(\Kb/K)$-modules are
 trivial, so this category won't do.
\end{remark}

{\bf Proof of Proposition 6}: It is sufficient to show that for
each $\ell$, the diagram of pairings:
$$\matrix{T^{-1}_j&\times&T^{1\vee}_{j-1}&\to&K^{*({\ell})}&\to&\bZl\cr
\|&&\|&&&&\uparrow\cr
T^{-1}_j&\times&T^{1\vee}_{j-1}&\to&&&\bZ.}$$ is commutative. The
pairing in the top row is that just defined, and the pairing in
the bottom row is constructed from the isogeny
$$\widetilde{N}:\, T^{-1}_j\to T^1_{j-1}$$
defined in Lemma 5. This shows that the top pairing takes values
in $\bZ$, hence factors through the image of $K^*$ in
$K^{*(\ell)}$ for every $\ell$. Clearly, we only need to prove
that this diagram is commutative after tensoring by $\bQ$.

Now, the commutativity of the diagram is deduced from Lemma 5 by
using the same argument as in the paper of Raynaud ([Ray2],
section 4.6). Concretely, given an exact sequence $0\to E' \to E
\to E/E' \to 0 $ as Galois $\bQl$-modules such that $E'(-1)$ and
$E/E'$ are $\bQl$-vector spaces of finite dimension with trivial
Galois action, one can construct as before a pairing
$$\{-,-\}:E/E' \times E'(-1)^{\vee} \to K^{*(\ell)}\otimes_{\bZl}\bQl
\stackrel{v_{\ell}}{\to}\bQl.$$ Then the map $N'\colon E/E' \to
E'(-1)$ constructed from this pairing coincides with the map
induced from the monodromy on $E$.\\

\begin{definition}If the Galois action on $T^{-1}_j(\Yb)$ and $T^1_{j-1}(\Yb)$ is trivial, we define the $j$-th intermediate Jacobian $J^j(X)$ as follows:
consider the pairing that we just defined:
$$\{-,-\}:\Lambda\times \Lambda^{'\vee}\to
K^*.$$  Now define:
$$J^j(X):=\mbox{Hom}(\Lambda^{'\vee}, {\bf G}_m) /\Lambda ,$$
where we are viewing $\Lambda$ as a subgroup of
$\mbox{Hom}(\Lambda^{'\vee}, {\bf G}_m)$ via the pairing $\{ -, -
\}$. This quotient makes sense in a rigid-analytic setting due to
Proposition 2, which tells us that $\Lambda$ is a lattice in the
algebraic torus $\mbox{Hom}(\Lambda^{'\vee}, {\bf
G}_m)$.\\
\end{definition}

Now we consider the case where there may be a non-trivial action
of the absolute Galois group $G_F$ on $\free{T^{-1}_{j}}$ and on
$\free{T^{1}_{j-1}}$, necessarily through a finite quotient. After
a finite unramified Galois extension $L/K$, we have a pairing
$$\free{T^{-1}_j}\times \free{T^{1\vee}_{j-1}}\to L^*.$$
If we show that this pairing respects the action of the Galois
group of $L/K$ on both sides, we will have an identification of
$\free{T^{-1}_{j}}$ and a lattice in the torus with character
group $\free{T^1_{j-1}}$. We can then define the $J^j(X)$ as the
\textit{non-split} analytic torus defined as the quotient. But
this pairing was defined by using an extension in
$\mbox{Ext}^1_G(\free{T^{-1}_{j}}\otimes\bZl,\free{T^{1}_{j-1}}\otimes\bZl(1)
)$ for any $\ell$, so it is clear that pairing we get
$$\{-,-\}:\free{T^{-1}_{j}}\times \free{ T^{1\vee}_{j-1}}\to \prod _{\ell}
L^{*(\ell)}$$ respects the Galois action of $\mbox{Gal}(L/K)$.

\begin{remark}
\begin{itemize}
\item[(i)]Observe that the dimension of $J^j(X)$ is equal to the rank of
$T^{1}_{j-1}$, which is equal to
$\dim_{\bQp}(H^{j}(\Yb,W\omega^{j-1}_{\Yb,log})\otimes_{\bZp}\bQp)$
(see [RX], Corollary 5).  This last dimension is equal to
$h^{j,j-1}:=\dim _{K}(H^{j}(X,\Omega^{j-1}_{X}))$ by ([I3],
Proposition 2.6 (c)).

\item[(ii)]Note that our intermediate Jacobian, which  \textit{a priori}
depends on the chosen regular proper model of $X$ over $R$,
actually only depends modulo isogeny only on the variety $X$. To
show this, observe first that the monodromy filtration we
constructed in the $p$-adic cohomology $H^n(\Xb,\bQp)$ is uniquely
defined, because there are no maps between $\bQp(j)$ and $\bQp(i)$
if $i\ne j$. On the other hand, the Tate module $V_p(J^j(X))$ is
isomorphic to $M_1(H)/M_{-3}(H)$, where $H=H^{2j-1}(\Xb,\bQp(j))$
and $M_{\bullet}$ is the monodromy filtration. But an analytic
torus over $K$ is determined modulo isogeny by its Tate module
(the proof of this is a generalization of the argument in [Se],
Ch. IV, \S A.1.2).  Observe that for $j=1$ and $j=d$, the
subquotient of $H^{2j-1}(\Xb,\bZl)$ that we use to construct the
intermediate Jacobian is, in fact, the whole torsion free quotient
of this cohomology group.  This fact, together with similar
arguments just given to determine an analytic torus by its Tate
modules,  then show that for $j=1$ and $j=d$, $J^j(X)$ is uniquely
determined by $X$ and is equal to, respectively, the Picard and
Albanese varieties of $X$.
\end{itemize}
\end{remark}

\section{The Abel-Jacobi mapping}

        Let $K$ be a finite extension of $\bQp$, and let $X$ be a smooth
projective variety over $K$.  We assume that $X$ has a regular
proper model $\cX$ over the ring of integers $R$ of $K$, with
special fibre satisfying the assumptions in \S 1.2.  In the last
section, we defined intermediate Jacobians $J^j(X)$, which we will
denote by $J^j$ when $X$ is fixed throughout the discussion.  In
this section, we define the Abel-Jacobi mapping:
$$CH^j(X)_{hom}\to J^j(X)\otimes_{\bZ}\bQ,$$
where $CH^j(X)_{hom}$ is the subgroup of the Chow group $CH^j(X)$
consisting of cycles in the kernel of the cycle map:
$$CH^j(X)\to \prod_{\ell} H^{2j}(\Xb,\bZl(j)),$$
where the product is taken over all prime numbers $\ell$.

        Let $T_{\ell}=T_{\ell}(J^j(X))$, the Tate module of $J^j(X)$,
and let $T=\prod_{\ell} T_{\ell}$.  Put
$V_{\ell}=T_{\ell}\otimes_{\bZl}\bQl$, and let $V=\prod_{\ell}
V_{\ell}$. Recall the groups $H^1_g(K,T)$, $H^1_g(K,V)$ of
Bloch-Kato ([BK2], 3.7).  We have
$$H^1_g(K,V_p)=\mbox{Ker}[H^1(K,V_p)\to
H^1(K,V_p\otimes_{\bQp}B_{DR})],$$ where $B_{DR}$ is the field of
``$p$-adic periods'' of Fontaine.  For $\ell\neq p$, we have:
$$H^1_g(K,V_{\ell})=H^1(K,V_{\ell}).$$
$H^1_g(K,T)$ is defined as the inverse image of $H^1_g(K,V)$ under
the natural map
$$H^1(K,T)\to H^1(K,V).$$

The following Lemma is proved as in the case of an abelian variety
in [BK2], \S 3.

\begin{lemma}
We have an isomorphism:
$$J^j(K)\to H^1_g(K,T),$$
compatible with passage to a finite extension of $K$ and with
norms.
\end{lemma}

\begin{remark}  The content of the lemma is that
the extensions of $\bZp$ by $\bZp(1)$ that give pairings with
values in $K^*$ (as opposed to $\widehat{K^*}$) are precisely the
de Rham extensions.
\end{remark}

Now, the following result is essentially due to Faltings [Fa2]: it
is a consequence of his results that the {\'e}tale cohomology of a
smooth variety (not necessarily projective) is de Rham in the
sense defined by Fontaine (see for example [I1]).

\begin{theorem}(Faltings) Let $X$ be a smooth projective variety.
Let $H:=\prod_{\ell} H^{2j-1}(\Xb,\bZl(j))$ and let
$cl':CH^j(X)_{hom} \to H^1(K, H)$ be the canonical {\'e}tale
Abel-Jacobi mapping (see [Bl] and [R1]). Then the image of $cl'$
is contained in $H^1_g(K,H)$.
\end{theorem}

\begin{corollary}
Let $X$ be a smooth projective variety over $K$ with a regular
proper model $\cX$ over R satisfying the assumptions of \S 1.2.
Let $M_{\cdot}$ be the monodromy filtration on
$H^{2j-1}(\Xb,\bQl(j))$, and let $V=\prod_{\ell} V_{\ell}$. Then
the $\ell$-adic Abel-Jacobi map $cl'$ factorizes through a map
$cl':CH^j(X)_{hom} \to H^1_g(K, M_1).$
\end{corollary}
{\bf Proof}: Observe first that for every $\ell$ we have that
$V_{\ell}/M_1$ is a successive extension of $\bQl$-vector spaces
of the form $\bQl(i)$, with $i<0$. Then the result follows from
the fact that $[V_{\ell}/M_1]^G=0$ and $H^1_g(K,\bQl(i))=0$ for
$i<0$ ([BK2], Example 3.5).\\

\begin{definition}
We define the rigid-analytic Abel-Jacobi mapping as the
composition
$$CH^j(X)_{hom}\to H^1_g(K,M_1)\to H^1_g(K,M_1/M_{-3})\cong J^j(X)(K)\otimes \bQ$$
 using Lemma 11 and Corollary 14.
\end{definition}

\begin{example}
Let $\cX$ be a regular proper scheme of relative dimension $d$
over $R$, whose special fibre satisfies the assumptions of \S 1.
Then $J^1(X)$ is the Jacobian variety of $X$, $J^d(X)$ is the
Albanese variety of $X$ and the maps defined above are the
classical Abel-Jacobi mappings, tensored with $\bQ$.
\end{example}

\begin{remark}
\begin{itemize}
\item[(i)]  Recall that by definition, $H^1_g(K,T)$ contains the torsion
subgroup of $H^1(K,T)$ for any finitely generated $\bZl$-module
$T$.  Now for $i<0$, $H^1(K,\bZl(i))$ has nontrivial torsion, in
general, and so it is not clear that we can define the Abel-Jacobi
mapping to $J^j(X)$ integrally.  Of course, we can define the map
on the inverse image of $H^1(K,T\cap M_1)$ under the Abel-Jacobi
mapping, but it is not clear how this subgroup changes as we
extend the base field $K$.

\item[(ii)]  As mentioned in Remark 10, the dimension of
$J^j(X)$ is the Hodge number $h^{j,j-1}$.  This differs in general
from the dimension of the complex intermediate Jacobian of
Griffiths, which is always equal to $\frac{B_{2j-1}}{2}$.  Using
higher odd K-theory, one can in some sense ``recover'' the
``rest'' of the classical intermediate Jacobian.  For example, for
$j=2$, we have the exact sequence:

$$0\to M_{-3}\to M_1\to M_1/M_{-3}\to 0.$$

The term on the right is the Tate module of our $J^2(X)$.  We then
get a natural map:

$$\ker \alpha_2\to H^1(K,M_{-3}),$$

where $\alpha_2$ is the Abel-Jacobi mapping we have just defined.
Now $M_{-3}$ modulo torsion is a direct sum of $\bZl(2)$'s, and we
have

$$\plim K_3(K,\bZ/\ell^n)\cong H^1(K,\bZl(2))$$

(see [L] and [MS]).

Thus we can relate $\ker\alpha_2$ to $K_3(K)$.  Using work of
Hesselholt-Madsen (see [HM], Theorem A), we can relate the kernel
of $\alpha_j$ successively to higher odd continuous algebraic
K-groups $\plim K_{2r-1}(K,\bZ/\ell^n),$ for $r=2,\cdots, j$. We
hope to pursue this point in another paper.
\end{itemize}
\end{remark}

\section{An Example}
        In this section we describe in detail an example that motivated
our theory.

\subsection{Product of Two Tate Elliptic Curves}
To warm up, we analyze the case of the product $X$ of two Tate
elliptic curves $E_1, E_2$, that is, elliptic curves with split
multiplicative reduction, or, equivalently, rigid analytic tori of
dimension 1. We can write $E_i=\bG_m/q_i^{\bZ}$, for some $q_i\in
K^*$ with $|q_i|<1$. The results we shall describe here are due to
Serre and Tate (see [Se], Ch. IV \S A.1.4), but we interpret them
in terms of our monodromy filtration, which will put them in a
form suitable for generalization to higher cohomology groups.\\

Let $V=H^1(\Eb_1,\bQp(1))$ and $W=H^1(\Eb_2,\bQp)$.  By the
K{\"u}nneth formula, we have that $H^2(\Xb,\bQp(1))$ is the direct
sum of $V\otimes W$ and two copies of the trivial representation
$\bQp$.  The monodromy filtration on $M=V\otimes W$ is the tensor
product of the monodromy filtrations on $V$ and $W$.  We thus have

$$M_0=\sum_{i+j=0}V_i\otimes W_j.$$

We have exact sequences:

$$0\to \bQp(1)\to V\to \bQp\to 0$$
$$0\to \bQp\to W\to \bQp(-1)\to 0,$$

which encode all of the information about the monodromy
filtrations. Let $q_i\, (i=1,2)$ be the class of the extension $V$
(resp. $W$) in
$$Ext^1_G(\bQp,\bQp(1))=K^{*\,(p)}\otimes_{\bZp}\bQp.$$

$q_i$ is in fact the image of the Tate parameter of $E_i$ under
the natural map:

$$K^*\to K^{*\,(p)}\otimes_{\bZp}\bQp,$$

which hopefully explains this abuse of notation.

Tensoring the first sequence with $W$ and the second by $V$, we
get exact sequences

$$0\to W(1)\to V\otimes W\to W\to 0.$$
$$0\to V\to V\otimes W\to V(-1)\to 0.$$

This gives us a map:

$$V\oplus W(1)\to V\otimes W,$$

and since $V=V_1, W=W_1$, we get a surjection

$$V\oplus W(1)\to M_0.$$

Note that the dimension of the space on the left is 4 and the
dimension of that on the right is 3.  Thus the kernel is of
dimension 1, and is easily seen to be $V\cap W(1)=\bQp(1)$, where
these two spaces are viewed as subspaces of $M_0$ via the exact
sequences above.  We also have that $M_{-2}=V_{-1}\otimes
W_{-1}=\bQp(1)$.  Summarizing, we have an exact sequence:

$$0\to\bQp(1)\to M_0\to M_0/M_{-2}\to 0,$$

where the Galois group acts on the right hand space through a
finite quotient.  Passing to an open subgroup $H$ of finite index
that acts trivially on this space and taking the standard basis of
$\bQp^2$, we see that the map

$$M_0/M_{-2}\to H^1(H,\bQp(1))$$

sends the vector $(a_1,a_2)$ to the class of $q_1^{a_1}q_2^{a_2}$
in $K^{*\,(p)}\otimes_{\bZp}\bQp$.  Thus $M_0^H\neq 0$ iff there
exist nonzero $p$-adic numbers $a_1, a_2$ such that
$q_1^{a_1}q_2^{a_2}=1$. We now have:

\begin{theorem}  (see [Se], Ch. IV, \S A.1.2) $E_1$ and $E_2$ are isogenous
iff there exist nonzero integers $a_1,a_2$ such that
$q_1^{a_1}q_2^{a_2}=1$.  The map

$$Hom(E_1,E_2)\otimes_{\bZ}\bQp\to Hom_G(V_p(E_1),V_p(E_2))$$

is an isomorphism.
\end{theorem}

Taking into account the other contributions to $H^2(\Xb,\bQp(1))$,
we get

\begin{corollary}
For $X$ the product of two Tate elliptic curves, the cycle map:

$$Pic(X)\otimes \bQp\to H^2(\Xb,\bQp(1))^G$$
is surjective.
\end{corollary}

This last statement is the ``$p$-adic Tate conjecture'' for the
product of two Tate
elliptic curves (see Remark 25(ii) below).\\

\subsection{Product of Three Tate Elliptic Curves}
Now we generalize this to $H^3$ of the product $X$ of
three Tate elliptic curves, $E_i, i=1,2,3$  with parameters $q_i$.
We study $J^2(X)$. Now $X$ is a 3-dimensional abelian variety, so
we have:

$$\mbox{dim}\,H^3(\Xb,\bQl)\cong\mbox{dim}\bigwedge^3H^1(X,\bQl)=\pmatrix{6\cr
3}=20.$$ Thus the complex intermediate Jacobian would be of
dimension 10.  Our $J^2(X)$ will be of dimension 9, and we explain
below how to recover the part that has been ``lost'' .\\

        We calculate the monodromy filtration on $H^3(\Xb,\bQp(2))$.
As above, the monodromy filtration on each
$V_i=H^1(\overline{E_i},\bQp(1))$ is given by:

$$V_{i,-1}=\bQp(1), V_{i,0}=\bQp(1), V_{i,1}=V,$$
and $V_i/V_{i,-1}=\bQp$.  The class of the extension

$$0\to V_{i,-1}\to V_i\to V_i/V_{i,-1}\to 0$$

in $H^1(K,\bQp(1))=K^{*\,(p)}\otimes_{\bZp}\bQp$ is given by the
class of the Tate parameter $q_i$ of $E_i$ in this group.  The
monodromy filtration on $M=V_1\otimes V_2\otimes V_3(-1)$ is the
tensor product of the monodromy filtrations on the $V_i$. Thus

$$M_n=\sum_{r+s+t=n}V_{1,r}\otimes V_{2,s}\otimes V_{3,t}(-1);$$
note that this sum is not direct.  We calculate easily the
following table, in which the top row is the level of the
filtration and the bottom row is the dimension:

$$\matrix{-3&-2&-1&0&1&2&3\cr
1&1&4&4&7&7&8}.$$

If we include the other contributions to $H^3(\Xb,\bQp(2))$ from
the K{\"u}nneth formula, we find the following table, which gives
the dimension of the spaces in the various steps of the monodromy
filtration on $H^3(\Xb,\bQp(2))$:

$$\matrix{-3&-2&-1&0&1&2&3\cr
1&1&10&10&19&19&20}.$$

The $p$-Tate vector space $V_{p}(J^2(X))$ of our $J^2(\Xb)$ is
$M_1/M_{-3}$, which is of dimension 18, and thus the dimension of
$J^2(X)$ is 9.  The ``lost'' dimension comes from the fact that
$Gr_{-2}M\cong \bQp(2)$ as Galois modules, and hence this cannot
contribute to the Tate module of a $p$-adic analytic torus.  As
mentioned above in Remark 17(ii), the lost part can be directly
related to $$K_3(K,\bQp)^{ind}=[\plim
K_3(K,\bZ/p^n)]\otimes\bQp,$$ via extensions of $\bQp$ by
$\bQp(2)$.\\

In some cases, we can also determine the image of the restriction
of the Abel-Jacobi mapping to $CH^2(X)_{alg}$; this is an abelian
variety, which is the universal abelian variety into which
$CH^2(X)_{alg}$ maps via a regular homomorphism (see [Mu]); we
denote it by $J^2_a(X)$. Its dimension should depend on the
multiplicative relations between the $q_i$ with integer exponents,
as we now explain.  To simplify the exposition,  {\bf we assume
that the Galois group acts trivially on the $T^i_j(\Yb)$}, where
$\Yb$ is the geometric special fibre of a strictly semi-stable
model $\cX$ of $X$ over the ring of integers of $K$. As we already
said, this will always happen after a finite unramified extension
of $K$ because all of the $T^i_j(\Yb)$ are already torsion free
abelian groups. \\

\subsection{Enriched monodromy operators}

With notation as in the previous sections, we have an extension
for each prime number $\ell$:

$$0\to T^3_0\otimes \bQl(1)\to M_{-1}\to T^1_1\otimes\bQl\to 0.$$

As in the definition of the intermediate Jacobian, we can define a
pairing:

$$(T^1_1\otimes_{\bZ}\bQ)\times (T^{3\vee}_0\otimes_{\bZ}\bQ)\to
K^{*(\ell)}\otimes_{\bZ}\bQ$$

by taking elements $\alpha\in T^1_1\otimes_{\bZ}\bQ, \beta\in
T^{3\vee}_0\otimes_{\bZ}\bQ$, tensoring the subgroups they
generate by $\bQl$ and $\bQl(1)$, respectively, pulling back the
extension above by the first and pushing it out by the second.
Taking the product for all $\ell$ and using Proposition 6 above,
we see that we get a $K^*\otimes_{\bZ}\bQ$-valued pairing and
hence a map:

$$N_1:\:  T^1_1\otimes_{\bZ}\bQ\to Hom(
T^{3\vee}_0\otimes_{\bZ}\bQ,K^*\otimes_{\bZ}\bQ),$$

that we call the {\it enriched monodromy operator}. See [R2] for
more on enriched monodromy operators and their relation with the
generalized Hodge-Tate conjecture.  We denote also by
$$N_1':\: T^{-1}_2\otimes_{\bZ}\bQ\to Hom(
T^{3\vee}_0\otimes_{\bZ}\bQ,K^*\otimes_{\bZ}\bQ),$$ the
composition of the map $N:\: T^{-1}_2\to T^1_1$ (tensored with
$\bQ$) with the map $N_1$ just defined.

\subsection{The Generalized Hodge-Tate Conjecture for $H^3$}Recall that if $Z$ and $X$ are two smooth projective
varieties over a field, then a correspondence between them is an
element $\eta$ of $CH^r(Z\times X)$;  such a cycle induces a
Galois equivariant map:

$$\eta_*:\:H^i(\Zb,\bQl(1))\to H^{i+2r}(\Xb,\bQl(1+r))$$

defined by $\eta_*(\alpha)=pr_{2*}(pr_1^*(\alpha)\cup [\eta]).$
Here $[\eta]$ is the cohomology class of $\eta$ in
$H^{2r}(\Zb\times \Xb,\bQl(r))$.  Let $N^rH^i(\Xb,\bQl(1+r))$ be
the subspace spanned by the images of such maps $\eta_*$, where
$Z$ is smooth and projective of dimension $d-r$, where
$d=\mbox{dim}(X)$.  This is called the {\it coniveau filtration}
on $X$.  Note the superscipts as we shall also use subscripts on
$N$

\begin{lemma}  Let $A$ be an abelian variety
contained in $J^2_a(X)$.   Write $V_l(A)$ as an extension

$$(*)\: 0\to S^1_0\otimes\bQl(1)\to V_l(A)\to S^{-1}_1\otimes\bQl\to
0,$$

where $S^i_j$ are finitely generated abelian groups (see the proof
below for why this is possible).  Then for some finite extension
$L/K$, we have:

$$S^1_0\otimes_{\bZ}\bQ\subseteq \ker N_{1,L}:\:
T^1_1\otimes_{\bZ}\bQ\to Hom(
T^{3\vee}_0\otimes_{\bZ}\bQ,L^*\otimes_{\bZ}\bQ)$$

and

$$S^{-1}_1\otimes_{\bZ}\bQ\subseteq \ker N'_{1,L}:\: T^{-1}_2\otimes_{\bZ}\bQ
\to Hom( T^{3\vee}_0\otimes_{\bZ}\bQ,L^*\otimes_{\bZ}\bQ).$$
\end{lemma}
{\bf Proof}: Since $A$ is contained in $J^2(X)$, it is an analytic
torus. Denote by $S^1_0$ and by $S^{-1}_1$ the corresponding
lattices of $A$, such that $A(L)\cong Hom(S^1_0,L^*)/S^{-1}_1$.
Now the inclusion map from $A$ to $J^2(X)$ gives us monomorphisms
$i\colon S^1_0 \to T^1_1$ and $i'\colon S^{-1}_1\to T^{-1}_2$. We
must show that $N_1\circ (i\otimes_{\bZ}\bQ)=0$ and $N'_1\circ
(i'\otimes_{\bZ}\bQ)=0$ for some finite extension $L/K$.  By the
definition of the enriched monodromy operator, to prove this it is
sufficient to show that for any prime $\ell $, the composition
$$S^1_0\otimes_{\bZ}\bQl \to T^1_1\otimes_{\bZ}\bQl
\stackrel{N_{L,\ell}}{\longrightarrow} Hom(
T^{3\vee}_0\otimes_{\bZ}\bQl,L^{*(\ell)}\otimes_{\bZl}\bQl)$$ is
trivial.

There exists a divisor $D$ on $X$ with desingularization $D'$ and
a correspondence $\eta\in CH^2(X\times D')$ such that
$V_{\ell}(A)$ is the image of $\eta_*$ in $H^3(\Xb,\bQl(2))$. Take
$L$ a finite extension of $K$ such that $D'$ has semi-stable
reduction. Then the cohomology $H^1(\overline{D'},\bQ\ell(1))$ is
semi-stable.

First consider the case $\ell \ne p$. Denoting by $M_{\bullet}$
the monodromy filtration on the cohomology, we have that the
correspondence $\eta$ gives a map respecting this filtration.
This gives us a map
$$Gr^M_{-1}(H^1(\overline{D'},\bQp(1))) \to S^{1}_{0}\otimes_{\bZ}\bQl \to
T^1_1\otimes_{\bZ}\bQl.$$
We have, on the other hand, a
commutative diagram of monodromy operators
$$\matrix{T^1_1\otimes_{\bZ}\bQl&\stackrel{N_{1,L}}{\to}&
Hom(
T^{3\vee}_0\otimes_{\bZ}\bQl,L^{*(\ell)}\otimes_{\bZl}\bQl)\cr
\uparrow&&\uparrow\cr Gr^M_{-1}(H^1(\overline{D'},\bQp(1)))
&\stackrel{N_{1,L}}{\to}&
Hom(Gr^M_{-3}(H^1(\overline{D'},\bQp(1)))^{\vee}(1),L^{*(\ell)}
\otimes_{\bZl }\bQl).}$$ But
$Gr^M_{-3}(H^1(\overline{D'},\bQp(1)))=0$, since there is no
graded piece of weight $\le -3$ in the monodromy filtration on
$H^1(\overline{D'},\bQp(1))$, and hence the image of
$S^1_{-1}\otimes_{\bZ}\bQl$ in $T^1_1\otimes_{\bZ}\bQl$ is
contained in $\ker N_L$, as claimed.

Now, to show the case $\ell=p$, we use instead of the $p$-adic
cohomology of $D'$ the log-crystalline cohomology of a suitable
model of $D'$. We have then, by applying Tsuji's theorem if
necessary, that the correspondence $\eta$ gives us a map between
the respective log-crystalline cohomologies, respecting the
monodromy filtration. Then the same argument applies as we have
just given for $\ell\neq p$. A similar argument applies also to
the
other claimed inclusion.\\

Conversely, we have the

\begin{conjecture}  ($p$-adic Generalized Hodge conjecture for
$H^3$)For $X$ with totally degenerate reduction, let the enriched
monodromy operators $N_{1,L}$ and $N'_{1,L}$ be defined as in \S
4.3. Then we have

$$Gr^{-1}N^1H^3(\Xb,\bQp)\!=\!$$
$$\sum_{[L:K]<\infty}\!\![\ker N_{1,L}:\: T^1_1\otimes_{\bZ}\bQ\to
Hom(
T^{3\vee}_0\otimes_{\bZ}\bQ,L^*\otimes_{\bZ}\bQ)]\otimes_{\bQ}
\bQp(-1).$$

and

$$Gr^M_{-1}N^1H^3(\Xb,\bQp)\!=\!$$
$$\sum_{[L:K]<\infty}\!\![\ker N'_{1,L}:\:T^{-1}_2\otimes_{\bZ}\bQ\to
Hom(
T^{3\vee}_0\otimes_{\bZ}\bQ,L^*\otimes_{\bZ}\bQ)]\otimes_{\bQ}
\bQp(-2).$$

Note that these spaces are of the same dimension, since the
monodromy operator $N:\:T^{-1}_2\to T^1_1$ is an isogeny. This is
as it should be, since $N^1H^3(\Xb,\bQp(2))$ is analogous to a
pure Hodge structure of odd weight, and so should be even
dimensional.
\end{conjecture}

\begin{remark}
To explain the motivation for phrasing this conjecture as we have,
consider the extension:

$$0\to M_{-3}\to M_{-1}\to M_{-1}/M_{-3}\to 0,$$

where $M=V_p(E_1)\otimes V_p(E_2)\otimes V_p(E_3)(-1)$.

Taking $H$-cohomology for an open subgroup $H$ of finite index in
$G$, we get an exact sequence:

$$0\to M_{-1}^H\to [M_{-1}/M_{-3}]^H\stackrel{\partial_H}{\to}
H^1(H,M_{-3}).$$

The space in the middle is isomorphic to $(T^1_1\otimes \bQp)^H$.
Arguing as in \S 4.1, we see that the space $M_{-1}/M_{-3}$ is
three dimensional, and for sufficiently small $H$, the map
$\partial_H$ takes a vector $(a_1,a_2,a_3)$ to $\prod q_i^{a_i}$,
viewed as an element of $L^{*\,(p)}$. Thus $\ker\,\partial_H$
corresponds to relations between the $q_i$ in $L^{*\,(p)}$
\textit{with $p$-adic exponents}.  There is no reason why these
have to be relations with rational integer exponents, as they are
for the case of two $q_i's$. For example, suppose we choose a
uniformizing parameter $\pi$ of $K$ and a logarithm with
$\log(\pi)=0$.  Then finding $p$-adic numbers $a_1, a_2, a_3$ with

$$\prod_{i=1}^3q_i^{a_i}=1$$

is equivalent to solving the simultaneous equations:

$$\sum_{i=1}^3a_i\log(q_i)=0$$

and

$$\sum_{i=1}^3a_iv_i=0$$

in $\bQp$, where $v_i=v(q_i)$.  This is a system of two equations
in three unknowns, so there should always be a solution other than
$(0,0,0)$.  But in general, there is no solution to these
equations with the $a_i$ rational integers. In other words, the
kernel of the enriched monodromy operator:

$$T^1_1\to Hom( T^{3\vee}_0,K^*)$$

tensored with $\bZp$, may be strictly smaller than the kernel of
the map:

$$T^1_1\otimes\bZp\to Hom(T^{3\vee}_0,K^{*\,(p)}).$$
Note that this map is the completion of the first map with respect
to subgroups of finite index, not the tensor product of this map
with $\bZp$. This is why its kernel can be bigger, since the
natural map

$$K^*\otimes_{\bZ}\bZp\to K^{*\,(p)}$$

is surjective but not injective.  In [R2], this difference between
relations with $p$-adic integer exponents and with rational
integer exponents is conjecturally explained via the difference
between two types of coniveau filtrations on \'etale cohomology.
One filtration takes $\bQp$-coefficients throughout and the other,
which is in general bigger, takes the inverse limit of the
coniveau filtrations with $\bZ/p^n\bZ$-coefficients and tensors
with $\bQp$.
\end{remark}

In the following, we shall prove Conjecture 21 in some cases when
$X$ is the product of three Tate elliptic curves.  Complex
analogues of these results are well-known (see e.g. [Li], p.
1197). The strategy is simple: by Lemma 20, for $H$ a sufficiently
small open subgroup of $G$ with fixed field $L$, the dimension of
$\ker N_{1,L}$ provides an upper bound for the dimension of
$Gr^{-1}N^1H^3(\Xb,\bQp)$. For $X$ as above, we can easily compute
this bound, and in some cases we can explicitly produce divisorial
correspondences $\eta$ such
that the sum of the images of $\eta_*$ is of the required dimension.\\

\begin{conjecture}  (Compare 7.5 on p. 1197 of [Li]) Let $X=E_1\times
E_2\times E_3$ be a product of three Tate elliptic curves with
parameters $q_i\, (i=1,2,3)$, and let $r$ be the rank of the space
of triples of integers $(n_1,n_2,n_3)$ with

$$\prod_{i=1}^3q_i^{n_i}=1.$$

Then

$$dim J_a^2(X)=6+r.$$

For example, for ``generic'' $q_i$ (no multiplicative relations),
the dimension is 6, and if all of the $q_i$ are equal, the
dimension is 8.  In particular, for $X$ the product of three Tate
elliptic curves, the restriction of the Abel-Jacobi map to
$CH^2(X)_{alg}$ is never surjective onto $J^2(X)$.

\end{conjecture}

\begin{proposition}  The conjecture is true if there are two independent
multiplicative relations between the $q_i$,  if there is one
relation and two of the three are isogenous, or if there are no
nontrivial multiplicative relations between the $q_i$.
\end{proposition}

 {\bf Proof}:   The dimension is at least 6 since $J^2(X)$ contains two
copies of $X$, via algebraic cohomology classes of the type

$$\alpha\otimes\beta\otimes\gamma\in H^0(\Eb_i)\otimes H^1(\Eb_j)\otimes
H^2(\Eb_k),$$

for each permutation $(i,j,k)$ of $(1,2,3)$.\\

   To determine the dimension of $J^2_a(X)$, assume first that the rank of
the space of relations is two.  Since the valuations of the $q_i$
are all positive, all of the 2 by 2 minors of the two vectors must
be nonzero. We can row reduce chosen basis vector for the relation
exponents so that each contains exactly one 0, and these are in
different coordinate positions. But then by Theorem 18 in \S 4.1,
we have two relations between two of the $q_i$'s, which gives
isogenies between the respective $E_i$. Then it is easy to see
that the dimension of $J^2_a(X)$ is 8, because $[M_{-1}/M_{-3}]^G$
is then spanned by three cohomology classes of the type:

$$\alpha_{ij}\otimes \beta\in H^1(\Eb_i)\otimes H^1(\Eb_j)\otimes
H^1(\Eb_k),$$

where $\alpha_{ij}$ is the class of the graph of a nonzero isogeny
between $E_i$ and $E_j$. These classes cannot all lift to
$M_{-1}^G$, as then the map:

$$[M_{-1}/M_{-3}]^G\to H^1(G,M_{-3})$$

would be zero, which can't be, since there are linear combinations
of
the valuations of the $q_i$ that are nonzero.\\

 If there is say a multiplicative relation between $q_1$ and $q_2$,
and these have no multiplicative relations with $q_3$, then we can
produce another cohomology class supported on the graph of the
isogeny, as
above.\\

  If there are no relations, then the dimension is at most six, and so is
exactly six.
  This completes the proof of the proposition.

\begin{remark}
\begin{itemize}
\item[i)] When there is just one multiplicative relation which
does not arise from an isogeny between two of the curves, there is
no obvious way to produce the abelian variety of the expected
dimension in $J^2(X)$, but Conjecture 21 predicts that it should
exist.

\item[ii)]  In a sequel to this paper [R2], we will formulate a $p$-adic
generalized Tate conjecture, which predicts the dimension of the
coniveau filtration $N^iH^j(\Xb,\bQp)$ for any $i, j$ in terms of
the kernels of enriched monodromy operators.

\item[iii)]  Since the image of the Abel-Jacobi mapping restricted to
$CH^2(X)_{alg}$ is never all of $J^2(X)$  when $X$ is the product
of three Tate elliptic curves, we hope to be able to use the
quotient of $J^2(X)$ by this image to detect codimension two
cycles that are homologically equivalent to zero, but not
algebraically equivalent to zero.  A regular proper model for such
products has been constructed by Gross-Schoen ([GS], \S 6,
especially Props. 6.1.1 and 6.3, and Cor. 6.4), but we do not
believe it is strictly semi-stable. Further blowing-up should make
it such, however.  See also the paper of Hartl [Ha], where he
constructs strictly semi-stable regular models of ramified base
changes of varieties with strictly semi-stable reduction, and of
products of such.
\end{itemize}
\end{remark}

\centerline{\bf References}

\begin{itemize}

\item[{\bf [Be]}]  A. Besser  A generalization of Coleman's
$p$-adic integration, Inventiones Math. 142 (2000) 397-434.

\item[{\bf [Bl]}]  S. Bloch   Algebraic cycles and values of
L-functions I, J. f. die reine und angewandte Math. 350 (1984)
94-108.

\item[{\bf [BK1]}] S. Bloch and K. Kato, $p$-adic {\'e}tale
cohomology, Publ. Math. I.H.E.S. 63 (1986) 107-152.

\item[{\bf [BK2]}] S. Bloch and K. Kato, L-functions and Tamagawa
numbers of motives, in {\it The Grothendieck Festschrift, Vol I},
Progr.Math. 87, Birkh{\"a}user, Boston, 1990.

\item[{\bf [BGS]}] S. Bloch, H. Gillet and C. Soul{\'e}, Algebraic
cycles on degenerate fibers, in Arithmetic Aspects of Algebraic
Geometry, Cortona 1994, F. Catanese editor, 45-69.


\item[{\bf [Con]}]  C. Consani, Double complexes and Euler
factors, Comp. Math. 111 (1998) 323-358.


\item[{\bf [DeJ]}]  A.J. de Jong, Smoothness, alterations and
semi-stability, Publ. Math. I.H.E.S. 83 (1996) 52-96.

\item[{\bf [De1]}]  P. Deligne, Th{\'e}orie de Hodge I, Actes du
Congr{\`e}s International des Math{\'e}maticiens, Nice, 1970, Tome
I, 425-430.


\item[{\bf [dS]}]  E. de Shalit, The $p$-adic monodromy-weight
conjecture for $p$-adically uniformized varieties,  Compos. Math.
141  (2005),  no. 1, 101-120.

\item[{\bf [Fa1]}]  G. Faltings, $p$-adic Hodge theory, Journal of
the American Mathematical Society, 1 (1988) 255-299.

\item[{\bf [Fa2]}]  G. Faltings, Crystalline cohomology and
$p$-adic Galois representations, in Algebraic Analysis, Geometry
and Number Theory, Baltimore MD 1988, 25-80, Johns Hopkins
University Press, 1989.


\item[{\bf [GN]}] F. Guill{\'e}n and V. Navarro-Aznar, Sur le
th{\'e}or{\`e}me local des cycles invariants, Duke Math. J. 61
(1990), 133-155.

\item[{\bf [Gr]}] P. Griffiths, On the periods of certain rational
integrals: I and II, Ann. Math. 90 (1969), 460-541.

\item[{\bf [GS]}] B. Gross and C. Schoen, The modified diagonal
cycle on a triple product of curves, Ann. de l'Institut Fourier 45
(1995) 649-679.

\item[{\bf [Gro]}]  A. Grothendieck, Hodge's general conjecture is
false for trivial reasons, Topology 8 (1969) 299-303.

\item[{\bf [Ha]}]  U. Hartl, Semi-stability and base change,
Archiv der Mathematik 77 (2001), 215-221.

\item[{\bf [HM]}]  L. Hesselholt and I. Madsen, On the K-theory of
local fields, Annals of Math. 158 (2003)  1-113.


\item[{\bf [I1]}] L. Illusie, Cohomologie de de Rham et
cohomologie {\'e}tale $p$-adique (D'apr{\`e}s G. Faltings, J.-M.
Fontaine et al.) in {\it S{\'e}minaire Bourbaki}, 1989-1990,
Asterisque 189-190 (1990), 325-374.

\item[{\bf [I2]}]  L. Illusie, Ordinarit{\'e} des intersections
compl{\`e}tes g{\'e}n{\'e}rales, in {\it The Grothendieck
Festschrift, Vol II}, Progr. Math, Vol. 87, Birkh{\"a}user,
Boston, 1990.

\item[{\bf [I3]}] L. Illusie, R{\'e}duction semi-stable ordinaire,
 cohomologie {\'e}tale $p$-adique et cohomologie de de Rham, d'apr{\`e}s
Bloch-Kato et Hyodo, in {\it P{\'e}riodes p-adiques, S{\'e}minaire
de Bures, 1988} (J.-M. Fontaine ed.), Ast{\'e}risque 223 (1994),
209-220.


\item[{\bf [It]}]  T. Ito, Weight-Monodromy conjecture for
$p$-adically uniformized varieties, Invent. Math.  159  (2005),
no. 3, 607--656.


\item[{\bf [Le]}]M. Levine The indecomposable $K_3$ of fields,
Annales Sci. \'Ecole Normale Sup\'erieure 22 (1989) 255-344.

\item[{\bf [Li]}]  D. Lieberman, On higher Picard varieties,
American J. of Mathematics, 90 (1968) 1165-1199.

\item[{\bf [Ma]}]  Y. Manin, Three-dimensional hyperbolic geometry
as $\infty$-adic Arakelov geometry, Invent. Math. 104 (1991)
223-244.

\item[{\bf [MD]}]  Y. Manin and V. Drinfeld, Periods of $p$-adic
Schottky groups, J. reine und ang. Math 262/263 (1973) 239-247.

\item[{\bf [MS]}] A.S. Merkur'ev and A. Suslin, The group $K_3$ of
a field, Izv. Akad. Nauk. SSSR (1989), English translation in
Mathematics of the USSR:  Izvestia (1989).


\item[{\bf [Mu1]}]  D. Mumford, An analytic construction of
degenerating curves over a complete local ring, Compositio Math.
24 (1972) 129-174.

\item[{\bf [Mu2]}]  D. Mumford, An analytic construction of
degenerating abelian varieties over complete local rings,
Compositio Math. 24 (1972) 239-272.

\item[{\bf [Mur]}]  J. Murre, Un r{\'e}sultat en th{\'e}orie de
cycles alg{\'e}briques de codimension deux, Comptes Rendus de
l'Acad{\'e}mie des Sciences 296 (1983) 981-984.


\item[{\bf [RZ1]}] M. Rapoport and Th. Zink, {\"U}ber die locale
Zetafunktion von Shimuravariet{\"a}ten, Monodromie-filtration und
verschwindende Zyklen in ungleicher Characteristic, Invent. Math.
68 (1982), 21-101.


\item[{\bf [R1]}]  W. Raskind, Higher $\ell$-adic Abel-Jacobi
mappings and filtrations on Chow groups, Duke Math. J. 78 (1995)
33-57.

\item[{\bf [R2]}]  W. Raskind, A generalized Hodge-Tate conjecture
for varieties with totally degenerate reduction over $p$-adic
fields,  in Algebra and Number Theory : Proceedings of the Silver
Jubilee Conference University of Hyderabad, edited by Rajat
Tandon. New Delhi, Hindustan Book Agency.

\item[{\bf [R3]}]  W. Raskind, On the $p$-adic Tate
conjecture for divisors on varieties with totally degenerate
reduction, in preparation.

\item[{\bf [RX]}]  W. Raskind and X. Xarles, On the \'etale
cohomology of algebraic varieties with totally degenerate
reduction over $p$-adic fields, to appear in J. of Mathematical
Sciences of the University of Tokyo.

\item[{\bf [Ray1]}] M. Raynaud, Vari{\'e}t{\'e}s ab{\'e}liennes et
g{\'e}om{\'e}trie rigide, Actes du ICM, Nice, 1970, Tome I,
473-477.

\item[{\bf [Ray2]}] M. Raynaud, 1-motifs et monodromie
g{\'e}om{\'e}trique, in {\it P{\'e}riodes p-adi\-ques,
S{\'e}minaire de Bures, 1988} (J.-M. Fontaine ed.),
Ast{\'e}risque 223 (1994), 295-319.

\item[{\bf [Se]}] J.-P. Serre, Abelian $\ell$-adic representations
and elliptic curves, Benjamin 1967.


\item[{\bf [St]}] J.H. Steenbrink, Limits of Hodge Structures,
Invent. Math. 31 (1976), 223-257.

\item[{\bf [Tsu]}] T. Tsuji, $p$-adic {\'e}tale cohomology and
crystalline cohomology in the semistable reduction case,
Inventiones Math 137 (1999) 233-411.

\item[{\bf [W1]}] A. Weil, Foundations of Algebraic Geometry,
American Mathematical Society Colloquium Publications Series,
Volume 29, Revised Edition, 1962.

\item[{\bf [W2]}]  A. Weil, On Picard varieties, American Journal
of Math., 74 (1952) 865-894; also in Oeuvres Scientifiques, Volume
II, Springer-Verlag Berlin 1978.

\end{itemize}
\textit{Authors' addresses}
\vspace{.1in}\\
\textsc{Wayne Raskind:  Dept. of Mathematics\\University of
Southern California\\Los
Angeles, CA 90089-1113, USA}\\
\texttt{email:  raskind@math.usc.edu}\\

\textsc{Xavier Xarles:  Departament de
Matem{\`a}tiques\\Universitat Aut{\`o}noma de
Barcelona\\08193 Bellaterra, Barcelona, Spain}\\
\texttt{email:  xarles@mat.uab.es}

\end{document}